\documentclass[12pt, a4paper]{article}

\usepackage{amsfonts,amssymb,amsmath,amsthm,textcomp,mathrsfs}
\usepackage{graphicx}
\usepackage{hyperref}
\usepackage{titlesec}
\usepackage{bbm}
\usepackage{breqn}
\titlelabel{\thetitle.\quad}
\titleformat{\section}[block]{\Large\bfseries\filcenter}{\thesection.}{0.4em}{}
\titleformat{\subsection}[hang]{\bfseries\filright}{\thesubsection}{1em}{}

\hypersetup{
    colorlinks,
    linkcolor={blue!50!black},
    citecolor={blue!50!black},
    urlcolor={blue!80!black}
}
\usepackage{longtable}
\usepackage{csquotes}
\usepackage{color,xcolor}
\definecolor{linkclr}{rgb}{0,0,0.8}
\definecolor{citeclr}{rgb}{0,0.5,0}
\definecolor{urlclr}{rgb}{0.8,0.4,0.4}
\definecolor{Plum}{HTML}{8E4585}
\definecolor{Mahogany}{HTML}{C04000}
\definecolor{Maroon}{HTML}{800000}
\definecolor{BrickRed}{HTML}{841F27}

\renewcommand{\thesection}{\textcolor{black}{\arabic{section}}}
\renewcommand{\thesubsection}{\textcolor{black}{\arabic{section}.\arabic{subsection}}}

\usepackage[PetersLenny]{fncychap} % girish chose this
\ChNameVar{\fontsize{18}{12}\usefont{OT1}{phv}{m}{n}\selectfont}
\ChNumVar{\fontsize{60}{62}\usefont{OT1}{ptm}{m}{n}\selectfont\hspace{2pt}}
\ChTitleVar{\centering\Huge\sf}\ChRuleWidth{0.5pt}\ChNameUpperCase
\usepackage{geometry}
\newtheorem{thm}{Theorem}[section]
\newtheorem{cor}[thm]{Corollary}
\newtheorem{lem}[thm]{Lemma}
\newtheorem{prop}[thm]{Proposition}

\theoremstyle{definition}

\addtocounter{section}{0}

 2

\usepackage{setspace}
\usepackage{appendix}

 \newcommand\blfootnote[1]{%
   \begingroup
   \renewcommand\thefootnote{}\footnote{#1}%   \addtocounter{footnote}{-1}%
   \endgroup
 }

\begin{document}
\doublespacing

\begin{center}
{\Large {\bf On monochromatic representation of sums of squares of primes}} \\
\vspace{1mm}
{\em Kummari Mallesham, Gyan Prakash and D. S Ramana }
\end{center}

\vspace{5mm} 
\begin{abstract}
\noindent 
When the sequences of squares of primes  is coloured with $K$ colours, where $K \geq 1$ is an integer, let $s(K)$ be the smallest integer such that each sufficiently large integer can be written as a sum of no more than $s(K)$ squares of primes, all of the same colour. We show that $s(K) \ll K \exp\left(\frac{(3\log 2 + {\rm o}(1))\log K}{\log \log K}\right)$ for $K \geq 2$. This improves on $s(K) \ll_{\epsilon} K^{2 +\epsilon}$, which is the best available upper bound  for $s(K)$.

\end{abstract}


%\subjclass[2010]{Primary 11N36; Secondary 11P99}
%\keywords{additive energy, monochromatic sums of primes, large sieve.}

\section{Introduction}
\label{intro}

\noindent
A subset $ \mathcal{D} $ of $ \mathbf{N} $ is said to be an asymptotic basis of finite order if there exists a positive integer $ m $ such that every sufficiently large integer can be written as a sum of at most $ m $ elements of $ \mathcal{D} $. The smallest $ m $ for which the above property holds is called the order of the asymptotic basis $ \mathcal{D} $. There are two classical examples of asymptotic bases of finite order. The first example is provided by Lagrange's four squares theorem, which tells us that the set of squares is an asymptotic basis of order four (in fact, it is a basis of order four). The second one is the set of primes. A classical 3-primes theorem of Vinogradov asserts that every sufficiently large odd integer can be written as a sum of three primes; it follows from this theorem that the set of primes is an asymptotic basis of order at most four. 

\vspace{2mm}
\noindent
Given an asymptotic basis $ \mathcal{D} $, one can ask the following ``chromatic" question. Let $ K \geq 1 $ be an integer and let $ \{  \mathcal{D}_{i} : 1 \leq i \leq K \} $ be any colouring  of 
$ \mathcal{D} $ in $ K $ colours. Let $ s_{\mathcal{D}} (K) $ denote the smallest integer $ n $ (if it exists) such that every sufficiently large integer is expressible as a sum of at most $ n $ elements of $ \mathcal{D} $, all of the same colour. Then what is the optimal bound for $ s_{\mathcal{D}} (K) $ in terms of $ K $?

\vspace{2mm}
\noindent
S{\'a}rk{\"o}zy \cite{sar1} was the first to ask this question in the case when $ \mathcal{D} $ is either the set of squares or the set of primes.
N. Hegyv{\'a}ry and F. Hennecart \cite{heg} attacked the problems of S{\'a}rk{\"o}zy using elementary methods. They obtained the bounds $ s_{\mathcal{D}} (K) \ll (K \log K)^{5} $ when $ \mathcal{D} $ is the set of squares and $ s_{\mathcal{D}} (K) \ll K^{3} $ when $ \mathcal{D} $ is the set of primes. D. S. Ramana and O. Ramar{\'e} \cite{rr} obtained 
the bound $ s_{\mathcal{D}} (K) \ll K \log log{2K} $ when $ \mathcal{D} $ is the set of primes, improving on the bound given by N. Hegyv{\'a}ry and F. Hennecart. Until recently, the best known bound on $ s_{\mathcal{D}} (K) $ when $ \mathcal{D} $ is the set of squares was due to P. Akhilesh and D. S. Ramana \cite{akhi}, namely $ s_{\mathcal{D}} (K) \ll_{\epsilon} K^{2+ \epsilon} $. Very recently Gyan Prakash, D. S. Ramana and O. Ramar{\'e} \cite{gsr} obtained the bound 
$ s_{\mathcal{D}}(K) \ll K \exp\left(\frac{(3\log 2 + {\rm o}(1))\log K}{\log \log K}\right) $   when $ \mathcal{D} $ is the set of squares, which is the best possible upto a constant.

\vspace{2mm}
\noindent
A classical theorem of Hua asserts that the set of squares of primes is an asymptotic basis of finite order. In this paper, we consider a 
``chromatic" version of Hua's theorem, raised by F. Hennecart, and following ideas in \cite{gsr}  obtain the following theorem.

\begin{thm}
\label{th1}
For any integer $K \geq 2$ we have $s(K) \leq K \exp\left(\frac{(3\log 2 + {\rm o}(1))\log K}{\log \log K}\right)$.
\end{thm}

\noindent
Here $o(1) \ll \frac{\log \log \log K}{\log \log K}$ for all large enough $K$. 
This improves on the bound $s(K) \ll_{\epsilon}  K^{2+\epsilon} $ given by Theorem 1.1, page 181 of Guohua Chen \cite{chen}.

\vspace{2mm}
\noindent
Our path to Theorem \ref{th1} passes through the theorem below, which we state with the help of following notation. For any subset $ S $ of the squares of primes in the interval $ (N,4N] $, we shall write 

\vspace{-2mm}
\begin{equation}
\label{addenergy}
E_{6} (S) = \sum_{\substack{ p_{1}^2 + p_{2}^2 + \ldots + p_{6}^2 = p_{7}^2 + p_{8}^2 + \ldots + p_{12}^2, \\ p_{i}^2 \in S,  1 \leq i \leq 12}} \log p_{1}  \log p_{2}\ldots \log p_{12}.
\end{equation}

\begin{thm} 
\label{main2}
Let $A \geq e^{e^2}$ be real number. Then for all sufficiently large integers $N$, depending only on $A$, and any subset $S$ of the squares of primes in the interval $(N, 4N]$ with $|S| \geq \frac{N^{\frac{1}{2}}}{A \log N}$  we have

\begin{equation}
\label{three11}
E_{6}(S) \ll \frac{|S|^{11} (\log{N})^{11} }{N^{\frac{1}{2}}} \exp\left(\frac{\left(3\log 2 + o(1)\right)\log A}{\log \log A}\right),
\end{equation}

\noindent
where $o(1) \ll \frac{\log \log \log A}{\log \log A}$.

\end{thm}

\vspace{2mm}
\noindent
We prove Theorem \ref{main2} in Section \ref{seccir}. Indeed, we reduce the problem of bounding $ E_{6}(S) $ to the finite problem using the circle method: On the minor arcs estimations we use known exponential sum estimates over prime squares. Dealing with the major arcs we use methods of Sam Chow \cite{chow} to get the asymptotics on major arcs and reduce to the finite problem. We deal with the finite problem in Section
\ref{loc} and resolve the finite problem using a slightly modification of Gyan Prakash, D. S. Ramana and O. Ramar{\'e} \cite{gsr}. Theorem \ref{th1} is deduced from Theorem \ref{main2} by means of a classical application of the Cauchy-Schwarz inequality followed by an appeal to the argument from \cite{heg}, involving the finite addition theorem of S{\'a}rk{\"o}zy \cite{sar2}. We give details of this deduction in Section \ref{fin}.

\vspace{2mm}
\noindent
Throughout this article we use $e(z)$ to denote $e^{2\pi i z}$, for any complex number~ $z$ and write $e_{p}(z)$ for $e^{\frac{2 \pi i z}{p}}$ when $p$ is a prime number. Further, all constants implied by the symbols $\ll$ and $\gg$ are absolute 
except when dependencies are indicated, either in words or by subscripts to these symbols. The Fourier transform of an integrable function $f$ on ${\bf R}$ is defined by $\widehat{f}(u) = \int_{\bf R} f(t)e(-ut) dt$. Finally, the notations $[a,b], (a,b]$ etc. will denote intervals in ${\bf Z}$, rather than ${\bf R}$, with end points $a$, $b$ unless otherwise specified.

\section{The Finite Problem}

\label{loc}

\noindent
The main result of this section is Theorem \ref{th4}, which is a slightly modification of the Theorem 2.1 of \cite{gsr}. For the sake of completeness of our arguments we will provide proof of Theorem \ref{th4}in full details. Before state our main theorem in this section, we introduce the terms involving in this theorem. We shall suppose that $A \geq e^{e^2}$ and let $U = \prod_{p \leq w} p$, where $w = A^{25}$. In addition,  we let ${\mathcal Z}$ be a set of primes in the interval $ (\sqrt{N}, 2 \sqrt{N}] $  with    
   
\begin{equation}
\label{10.1a}
|{\mathcal Z}|  \geq \frac{\sqrt{N}}{A \log N} \;\;\text{and}\;\; 
|\{ z \in {\mathcal Z}| z \equiv a \,{\rm mod}\, U\}| \leq \frac{3\sqrt{N}}{\phi(U) \log N},
\end{equation}
   
\noindent
for all classes $a$ in ${\bf Z}/U{\bf Z}$. Also, we denote by ${\bf c} = \{c(i)\}_{i \in I}$ a given finite sequence of integers and let ${R}_{U}({\mathcal Z}, {\bf c})$ denote the set of triples $(x,y, i)$ in ${\mathcal Z} \times {\mathcal Z} \times I$ such that $x ^2 + y^2  + c(i)$ is an invertible square modulo $U$. Finally, let $\tau(U) = 2^{\pi(w)}$ be the number of divisors of $U$. Now we can state the theorem as follows.  

\begin{thm} 
\label{th4}
We have 
\vspace{-2mm}
\begin{equation}
\label{one2}
|{R}_{U}({\mathcal Z}, {\bf c})| \leq   \left(\frac{U}{\phi(U)} \right)^{2} \frac{|{\mathcal Z}|^2|I|}{\tau(U)} \exp\left(\frac{\left(3\log 2 + o(1)\right)\log 3A}{\log \log  3A}\right),
\end{equation}

\noindent
where $o \ll \frac{(\log \log \log 3A)}{\log \log 3A}$. 

\end{thm}

\noindent
We prove Theorem \ref{th4} in Subsection \ref{prlp}. The proof begins by using the optimisation principle given by Lemma \ref{bil} to pass to a problem in $\left( {\bf Z}/U{\bf Z} \right)^{*}$, considered in Theorem \ref{th3}. Two applications H{\"o}lder's inequality and the Chinese Remainder Theorem reduce the proof of Theorem \ref{th3} to the solution of a problem in $\left( {\bf Z}/p{\bf Z} \right)^{*}$ for a given prime $p|U$. This problem is treated by the Corollary \ref{modp} of the following  subsection.

\subsection{A  Sum in  {\bf Z}/p{\bf Z}}\label{subgauss}

\noindent
We write $G_p$ for the ring ${\bf Z}/p{\bf Z}$ when $p$ is a prime number. Also,  $\lambda_{p}(x)$ shall denote the Legendre symbol $(\frac{x}{p})$, for any $x$ in $G_p$.

\begin{lem}
\label{sumzp}
Let $p$ be a prime number and $c$ an element of $G_p$. Then for any an even integer $t \geq 2$ we have

\vspace{-2mm}
\begin{equation}
\label{02}
\sum_{{(y_1, y_2, \ldots,y_{t/2})} \in G_p^{t/2}} \sum_{(x_1, x_2, \ldots x_t) \in G_p^t} \prod_{\substack{1 \leq i \leq t,\\ 1\leq j \leq t/2.}} \left(1 + \lambda_{p}(x_i^2 + y_j^2 +c) \right) \leq 
p^{3t/2}\exp\left(\frac{4t^5 2^{t}}{p}\right).
\end{equation}

\end{lem}

\vspace{2mm}
\noindent
{\sc Proof.---} See the proof of Proposition 2.2 of \cite{gsr}, for example.

\vspace{2mm}
\noindent
Before stating a corollary of the above lemma, let us introduce some additional notation. Let $p$ be a fixed prime number and let $c$ be a given element of $G_p$. For any $(x,y)$ in $G_p^2$ we set $\delta_{p}(x,y) = \lambda_{p}(x^2+y^2 +c)$ and $\epsilon_{p}(x,y) = 1 +\delta_{p}(x,y)$. We endow $G_p \setminus \{ 0 \}$, and likewise  $ \left( G_p \setminus \{ 0 \} \right)^t$ for any integer $t\geq 1$, with their uniform probability measures and write ${\mathbb E}_{{x}}$ and  ${\mathbb E}_{x_1, x_2, \ldots, x_t}$ respectively in place of $\frac{1}{p-1}\sum_{x \in  G_p \setminus \{ 0\} }  $ and  $\frac{1}{(p-1)^t}\sum_{x_1, x_2, \ldots, x_t \in G_p \setminus \{ 0\}}$. Finally, we define ${\mathcal E}_{p}(k,t)$ for any integer $k$ with $1 \leq k \leq t$ by  

\vspace{-2mm}  
\begin{equation}
\label{0031e}
{\mathcal E}_{p}(k,t) = {\mathbb E}_{y_1, y_2, \ldots, y_t} {\mathbb E}_{x_1, x_2, \ldots, x_t} \prod_{\substack{1 \leq i \leq t,\\ 1\leq j \leq k.}} \epsilon_{p}(x_i,y_j) .
\end{equation}

\noindent
Using this notation we state our corollary as follows.

\begin{cor}
\label{modp}
For any even integer $ t \geq 2 $ we have 
\vspace{-2mm}  
\begin{equation}
\label{0031}
{\mathcal E}_{p}(t/2,t) \leq \left( \frac{p}{p-1}\right)^{2t} \exp\left(\frac{4t^5 2^{t} }{p}\right)\; .
\end{equation} 

\end{cor}

\vspace{2mm}
\noindent
{\sc Proof.---} Since $ t \geq 2 $ is an even integer. By taking $ k = t/2 $ in \eqref{0031e}  we get that

\vspace{-2mm}
\noindent

\begin{equation}
{\mathcal E}_{p}(t/2,t) = {\mathbb E}_{y_1, y_2, \ldots, y_t} {\mathbb E}_{x_1, x_2, \ldots, x_t} \prod_{\substack{1 \leq i \leq t,\\ 1\leq j \leq t/2.}} \epsilon_{p}(x_i,y_j).
\end{equation}

\noindent
Observing that the summands in the sum are independent of the variables $ y_{i} $ for $ t/2 < i \leq t $, and allowing the sum over full $ G_{p} $ we get

\begin{equation}
\label{cor1}
{\mathcal E}_{p}(t/2,t) \leq \frac{p^{t/2}}{(p-1)^{2t}} \sum_{{(y_1, y_2, \ldots,y_{t/2})} \in G_p^{t/2}} \sum_{(x_1, x_2, \ldots x_t) \in G_p^t} \prod_{\substack{1 \leq i \leq t,\\ 1\leq j \leq t/2.}} \left(1 + \lambda_{p}(x_i^2 + y_j^2 +c) \right). 
\end{equation}

\noindent
Substituting the upper bound on the sum in the right of \eqref{cor1} by Lemma \ref{modp}, we conclude the inequality \eqref{0031} holds.

\subsection{The Problem Modulo U} \label{modu}

\noindent
Let, as above, $ A \geq e^{e^2}$ be real number and  $U = \prod_{p \leq w} p$, where $w = A^{25}$. Suppose further that ${\mathcal X}$ and ${\mathcal Y}$ are subsets of  
$\left( {\bf Z}/U{\bf Z} \right)^{*} $ of density at least 
$\frac{1}{A}$. That is,  
   
\begin{equation}
\label{10.1a1}
|{\mathcal X}| \,\text{and}\, |{\mathcal Y}|  \geq \frac{\phi(U)}{A} \;.
\end{equation}
   
\noindent
For a given element $c$ of ${\bf Z}/U{\bf Z}$, let ${T}_{c}({\mathcal X},{\mathcal Y})$ denote the set of pairs $(x,y) \in {\mathcal X} \times {\mathcal Y}$ such that $x ^2 + y^2  + c $ is an invertible square in  ${\bf Z}/U{\bf Z}$.

\begin{thm} 
\label{th3}
For all $A$,  $U$, ${\mathcal X}$, ${\mathcal Y}$ and $c$ as above, we have 

\vspace{-2mm}
\begin{equation}
\label{one}
|{T}_{c}({\mathcal X},{\mathcal Y} )| \leq  \left(\frac{U}{\phi(U)} \right)^{2} \frac{|{\mathcal X}||{\mathcal Y}|}{\tau(U)} \exp\left(\frac{\left(3\log 2 + O\left(\frac{\log \log \log A}{\log \log A}\right)\right)\log A}{\log \log  A}\right).
\end{equation}

\end{thm}

\vspace{2mm}
\noindent
{\sc Proof.---}  We shall write $G$ for the set $ \left({\bf Z}/U{\bf Z} \right)^{*}$ and continue to use $G_p$ for ${\bf Z}/{p{\bf Z}}$. Also, for any $x$ in ${\bf Z}/U{\bf Z}$ and $p|U$ we denote the canonical image of $x$ in ${\bf Z}/p{\bf Z}$ by $x_p$ and, to be consistent with the notation of preceding subsection, write $\lambda_{p}(x)$ for the Legendre symbol $(\frac{x_p}{p})$. Then we have that

\vspace{-2mm}
\begin{equation}
\label{three}
|{T}_{c}({\mathcal X},{\mathcal Y})| \leq \sum_{x \in {\mathcal X}}\, \sum_{y \in \mathcal{Y}}  \prod_{p|U} \left(\frac{1 + \lambda_{p}(x^2+y^2+c)}{2}\right) \; , 
\end{equation}

\noindent
since $0 \leq 1 + \lambda_{p}(x^2+y^2+c) \leq 2$ for any pair $(x,y)$ in ${\mathcal X} \times {\mathcal Y}$, with equality in the upper bound for every prime $p|U$ when  $x^2+y^2+c$ is an invertible square in ${\bf Z}/U{\bf Z}$. On extending the definitions of $\delta_{p}$ and $\epsilon_p$ from Subsection \ref{subgauss} by setting  $\delta_{p}(x,y) = \lambda_{p}(x^2+y^2 +c)$ and $\epsilon_{p}(x,y) = 1 +\delta_{p}(x,y)$ for any $(x,y)$ in $ \left({\bf Z}/U{\bf Z} \right)^2$ and $p|U$, we may rewrite \eqref{three} as 

 \vspace{-2mm}
 \begin{equation}
 \label{three1}
 |{T}_{c}({\mathcal X},{\mathcal Y})| \leq \frac{1}{\tau(U)} \sum_{x \in {\mathcal X}}\, \sum_{y \in \mathcal{Y}}  \prod_{p|U} \epsilon_{p}(x,y) \; . 
 \end{equation}

\noindent
Let $t \geq 2$ be an even integer. Then an interchange of summations followed by an application of H\"older's inequality to exponent $t$ to the right hand side of \eqref{three1} gives

\vspace{-2mm}
\begin{equation}
\label{four}
|{T}_{c}({\mathcal X},{\mathcal Y})| \leq \frac{|{\mathcal Y}|^{1-\frac{1}{t}}}{\tau(U)} \left(\sum_{y \in {\mathcal Y}}\, \left(\sum_{x \in \mathcal{X}} \, \prod_{p|U} \epsilon_{p}(x,y)\right)^{t} \right)^{\frac{1}{t}} \;. 
\end{equation}

\noindent
To bound the sum over $y \in  {\mathcal Y}$ on the right hand side of the inequality above, we first expand the summand in this sum and extend the summation to all $y \in G$. By this we see that

\vspace{-2mm}
\begin{equation}
\label{five}
\sum_{y \in {\mathcal Y}}\, \left(\sum_{x \in \mathcal{X}} \, \prod_{p|U} \epsilon_{p}(x,y)\right)^{t} \leq 
\sum_{y \in G}\, \sum_{(x_1,x_2, \ldots,x_t) \in {\mathcal X}^{t}} \, \prod_{1 \leq i \leq t} \prod_{p|U} \epsilon_{p}(x_i, y) \, . 
\end{equation}

\noindent
Interchanging  the summations over $G$ and ${\mathcal X}^t$ on the right hand side of the above relation and applying H\"older's inequality again, this time to exponent $\frac{t}{2}$, we obtain that the right hand side of \eqref{five} does not exceed 

\vspace{-2mm}
\begin{equation}
\label{6}
 |{\mathcal{X}}|^{t-2} \left(\sum_{(x_1,x_2, \ldots,x_t) \in {\mathcal X}^{t}} \left(
\sum_{y \in G}\,  \prod_{1 \leq i \leq t} \prod_{p|U} \epsilon_{p}(x_i,y) \right) 
^{\frac{t}{2}}\right)^{\frac{2}{t}}. 
\end{equation}

\noindent
Finally, on expanding the summand in the sum over ${\mathcal X}^t$ in \eqref{6} and   extending the summation to all of $G^t$ we conclude using \eqref{five} and \eqref{four} and a rearrangement of terms that  

\begin{equation}
\label{7}
|{T}_{c}({\mathcal X},{\mathcal Y})| \leq \frac{|{\mathcal X}| |{\mathcal Y}|}{\tau(U)} \left(\frac{\phi(U)^3}{|{\mathcal X}|^2 |{\mathcal Y}|}\right)^{\frac{1}{t}} {\mathcal E}\left(\frac{t}{2},t\right)^{\frac{2}{t^2}}, 
\end{equation}

\noindent
where for any integer $k$ with $1 \leq k \leq t$ we have set  

\vspace{-2mm}
\begin{equation}
{\mathcal E}(k,t) = 
\frac{1}{\phi(U)^{2t}}\sum_{(y_1,y_2, \ldots,y_t) \in G^{t}} 
\sum_{(x_1,x_2,, \ldots, x_t) \in G^{t}}\, \prod_{p|U} \prod_{\substack{1 \leq i \leq t, \\ 1 \leq j \leq k.}}  \epsilon_{p}(x_i,y_j)\, .\end{equation}

\noindent
The Chinese Remainder Theorem gives $G = \prod_{p|U} \left( G_p \setminus \{ 0\} \right)$. Moreover, for all $p|U$ and $(x,y)$ in 
$ \left({\bf Z}/U{\bf Z} \right)^2 $  we have  $\epsilon_{p}(x,y) = \epsilon_{p}(x_p, y_p)$. It follows that ${\mathcal E}(k,t) = \prod_{p|U} {\mathcal E}_{p}(k,t)$, where ${\mathcal E}_{p}(k,t)$ is as defined by \eqref{0031e}. Using \eqref{0031} with $k = \frac{t}{2}$, valid on account of Corllary \ref{modp}, and recalling that $U = \prod_{ p \leq A^{25}} p$ we then obtain 

\vspace{-2mm} 
\begin{equation}
\label{8.5}
\left( {\mathcal E}(k,t) \right)^{\frac{2}{t^2}} = \left( \prod_{p|U} {\mathcal E}_{p}(k,t) \right)^{\frac{2}{t^2}} \leq \left(\frac{U}{\phi(U)} \right)^{4/t}\exp \left( 8 \,t^3\, 2^t \sum_{p \leq A^{25}} \frac{1}{p} \right).
\end{equation}

\noindent
From (3.20) on page 70 of \cite{rs} we deduce that $\sum_{p \leq A^{\ell}} \frac{1}{p}  \leq  
(\log {50}) \log\log A$, since $A \geq 4$. On combining this remark with \eqref{8.5}, \eqref{10.1a1} and \eqref{7} we then conclude that for any even integer $t \geq 2$ we have  

\begin{equation}
\label{8}
|{T}_{c}({\mathcal X},{\mathcal Y})| \leq \left( \frac{U}{\phi(U)} \right)^2 \frac{|{\mathcal X}| |{\mathcal Y}|}{\tau(U)}\exp\left( \frac{3\log A}{t} +  8\,(\log {50})\, t^3\, 2^t \log\log A \right). 
  \end{equation}

\noindent
Let us now set $v \log 2 = \log\left( \frac{\log A}{ (\log \log A)^6}\right)$ and suppose that $A_0 \geq e^e$ is such that we have $\frac{\log \log A}{\log \log\log A} \geq 12$ and $v \geq 4$ for all $A > A_0$. For such $A$ we take $t$ in \eqref{8} to be an even integer satisfying $v \leq t \leq v+2$. Also,  with $w = \frac{6\log\log\log A}{\log \log A}$ we have $w \leq \frac{1}{2}$  and $v = \frac{(1-w)\log\log A}{\log 2}$. Thus   $\frac{1}{t}  \leq \frac{1}{v} \leq \frac{(\log 2) (1+2w)}{\log\log A}$ and $t^32^t \leq 32v^3 2^v \leq \frac{32\log A}{(\log 2)^3(\log\log A)^3}$. Substituting these inequalities in \eqref{8} we obtain \eqref{one} for $A > A_0$. To obtain \eqref{one} for $e^{e^2} \leq A \leq A_{0}$ it suffices to take $t =2$ in \eqref{8}.

%For $ 4 \leq A \leq A_0$ we take $t=2$ and set $C(l)$ to be the supremum of the exponential factor on the right hand side of \eqref{8} taken over the set $t=2$ and $4 \leq A \leq A_0$. Then \eqref{8} tells us that \eqref{one} certainly holds for $ 4 \leq A \leq A_0$. When $A > A_0$ 
\qed

\subsection{ An Optimisation Principle} \label{opt}

\vspace{2mm}
\noindent
This subsection summarizes Subsection 2.3 of \cite{rr}.
Suppose that $n \geq 1$ is an integer 
and let $P$ and $D$  be real numbers $> 0$. Further, assume that 
the subset ${\mathcal K}$ of points $x = (x_1,x_2,\ldots,x_n)$ in ${\bf R}^n$ satisfying the conditions 

\begin{equation}
\label{opt1}
\sum_{1 \leq i \leq n} x_i = P \;\;\text{and}\;\; 0 \leq x_i \leq D \;\;\text{for all $i$ .}
\end{equation}

\noindent
is not empty. Then ${\mathcal K}$ is a non-empty, compact and  convex subset of ${\bf R}^n$ and we have the following  fact.

\begin{lem} \label{bil}
If $f: {\bf R}^n \times {\bf R}^n \mapsto {\bf R}$ a bilinear form with real coefficients $\alpha_{ij}$ defined by  $f(x,y) =\sum_{1\leq i, j \leq n} \alpha_{ij} x_i y_j$ then 

\vspace{1mm}
\noindent
$(i)$ There are extreme 
points $x^{*}$ and $y^{*}$ of ${\mathcal K}$ so that $f(x,x) \leq 
f(x^{*},y^{*})$ for all $x \in {\mathcal K}$.

\vspace{1mm}
\noindent
$(ii)$ If $x^{*} = (x_1^{*}, x_2^{*}, \ldots, x_n^{*})$ is an extreme point 
of ${\mathcal K}$ then, excepting at most one ~$i$, we have either $x_i^{*} = 0$ or $x_i^{*} = D$ for each $i$. Also, if $m$ is the 
number of $i$ such $x_i^{*} \neq 0$ then $mD \geq P > (m-1) D$.
 
\end{lem}

\vspace{2mm}
\noindent
{\sc Proof.---} See the proof of Proposition 2.2 of \cite{rr}, for example.

\subsection{Proof of Theorem \ref{th4}}\label{prlp}

 Let $a,b$ be any elements of ${\bf Z}/U{\bf Z}$. For any  $i$ in   $I$ we set $\alpha_i (a,b) = 1$ if $a^2 +b^2 +c(i)$ is an invertible square in  ${\bf Z}/U{\bf Z}$ and 0 otherwise. Further, we write $m(a)$ for the number of $z$ in  ${\mathcal Z}$ such that $z \equiv a$ mod $U$. Then if $\tilde{{\mathcal Z}}$ denotes the image of ${\mathcal Z}$ in ${\bf Z}/U{\bf Z}$ we have  

\begin{equation}
\label{th124}
|{R}_{U}({\mathcal Z}, {\bf c})|\, = 
 \, \sum_{i \in I} \sum_{(a,b) \in {\tilde{Z}}^2} \alpha_{i}(a,b)\, m(a)m(b) \; . 
\end{equation}

\noindent
Moreover, on account of the second assumption in \eqref{10.1a} we have that  

\begin{equation}
\label{th123}
\sum_{a \in \tilde{{\mathcal Z}}} m(a) = |{\mathcal Z}| \;\;\text{and}\;\; 0 \leq m(a) \leq D , 
\end{equation}

\noindent
where $ D =  \frac{3 \sqrt{N}}{\phi(U) \log N}$. For large values of 
$N$, depending on $A$;  $\tilde{{\mathcal Z}}$ is contained in 
 $ \left( {\bf Z}/U{\bf Z} \right)^{*}$. Let us bound the inner sum on the right hand side of \eqref{th124} for a fixed $i$ in $I$.  By means of Lemma \ref{bil} and \eqref{th123} we obtain   

\begin{equation}
\label{th125}
\sum_{(a,b) \in {\tilde{{\mathcal Z}}}^2} \alpha_{i}(a,b)\, m(a)m(b) \leq  \sum_{(a,b) \in {\tilde{Z}}^2} \alpha_{i}(a,b)\, x^{*}_{a} y^{*}_{b}
 \, ,
\end{equation}

\noindent
for some $x^{*}_{a}$ and $y^{*}_{b}$, with $a$ and $b$ varying over $\tilde{{\mathcal Z}}$, satisfying the following conditions. All
the $x^{*}_{a}$, and similarly all the $y^{*}_{b}$, are either 0 or $D$ excepting at most one, 
which must lie in $(0,D)$. Moreover, if ${\mathcal X}$ and ${\mathcal Y}$ are, 
respectively, the subsets of $\tilde{{\mathcal Z}}$ for which $x^{*}_{a} \neq 0$ and $y^{*}_{b} \neq 0$ 
then $|{\mathcal X}| D \geq |{\mathcal Z}| > (|{\mathcal X}|-1) D$. From the first condition in \eqref{10.1a} we then  get $|{\mathcal X}| \geq \frac{|{\mathcal Z}|}{D} \geq \frac{\phi(U)}{3A}\geq 2$. Consequently, we also have $D \leq \frac{|{\mathcal Z}|}{|{\mathcal X}|-1}\leq \frac{2|{\mathcal Z}|}{|{\mathcal X}|}$. The same inequalities hold with $|{\mathcal X}|$ replaced by $|{\mathcal Y}|$. Then with $T_{c(i)}({\mathcal X}, {\mathcal Y})$ as in  Subsection \ref{modu} we have that $\sum_{(a,b) \in {\mathcal X} \times {\mathcal Y}} \alpha_{i}(a,b) =|T_{c(i)}({\mathcal X}, {\mathcal Y})|$ and therefore that the right hand side of \eqref{th125} does not exceed 
$ \frac{4|T_{c(i)}({\mathcal X}, {\mathcal Y})||{\mathcal Z}|^2}
{|{\mathcal X}||{\mathcal Y}|} $. Using this in \eqref{th124} together with the bound supplied by \eqref{one} for $|T_{c(i)}({\mathcal X}, {\mathcal Y})|$, applicable since $3A \geq e^{e^2}$, we then conclude that \eqref{one2} holds.
 
\qed

\section{An Application of the Circle Method}
\label{seccir}

\noindent
We prove Theorem \ref{main2} in this section. As stated in Section \ref{intro}, we will first reduce the problem of bounding $ E_{6} (S) $ to the finite problem. This is carried out in Subsections \ref{mincon} through \ref{majcon} starting with the preliminaries given below. We then complete the proof of Theorem \ref{main2}  in Subsection \ref{maincom} by applying Theorem \ref{th4}. 

\vspace{2mm}
\noindent
We suppose that $A \geq e^{e^2}$ are real number and assume that $N$ is a sufficiently large integer depending only on $A$, its actual size varying to suit our requirements at various stages of the argument. We set $\alpha(t) = 1 -\left|\frac{2t}{5N}\right|$ when $|t|\leq \frac{5N}{2}$ and 0 for all other $t \in {\bf R}$ and set $\beta(t) = \alpha(t - \frac{5N}{2})$. Thus $\beta(t) \geq 0$ for all $t$ in ${\bf R}$ and $\beta(t) \geq \frac{2}{5}$ when $t \in [N, 4N]$. Finally, we set 

\vspace{-2mm}
\begin{equation}
\label{31}
\psi(t) = \sum_{n} \mathbbm{1}_{\mathbb{P}}(n) \,  2 \, n \,  \log{n} \, \beta(n^2) \, e(n^2t) \; 
\end{equation}

\noindent
and write $\widehat{S}(t) = \sum_{p^2 \in S} \log{p} \,  e(p^2t)$ for any $t \in {\bf R}$ for a given subset $S$ of the squares in $(N, 4N]$ satisfying the hypotheses of Theorem \ref{main2}. We observe that 

\vspace{-2mm}
\begin{equation}
\label{32}
\frac{4}{5} \sqrt{N}E_{6}(S) \, \leq\, \int_{0}^{1}  \widehat{S}(t)^6 \widehat{S}(-t)^{5}{\psi}(-t)\, dt \;.
\end{equation} 

\noindent
Indeed,  

\begin{align}
E_{6} (S) &= \sum_{\substack{ p_{1}^2 + p_{2}^2 + \ldots + p_{6}^2 = p_{7}^2 + p_{8}^2 + \ldots + p_{12}^2, \\ p_{i}^2 \in S.}} \log p_{1}  \log p_{2}\ldots \log p_{12}   \\
& \leq \sum_{\substack{ p_{1}^2 + p_{2}^2 + \ldots + p_{6}^2 = p_{7}^2 + p_{8}^2 + \ldots + p_{11}^2+q^2, \\ p_{i}^2 \in S; \, q^2 \text{ is a prime square in}  (N,4N].}} \log p_{1}  \log p_{2}\ldots \log p_{11} \log q  \label{e1} \\
& \leq 1/2 \sqrt{N} \sum_{\substack{ p_{1}^2 + p_{2}^2 + \ldots + p_{6}^2 = p_{7}^2 + p_{8}^2 + \ldots + p_{11}^2+q^2, \\ p_{i}^2 \in S; \, q^2 \text{ is a prime square in}  (N,4N].}} \log p_{1}  \log p_{2}\ldots \log p_{11} \, 2q \,  \log q  \label{e2} \\
& \leq 5/4 \sqrt{N} \sum_{\substack{ p_{1}^2 + p_{2}^2 + \ldots + p_{6}^2 = p_{7}^2 + p_{8}^2 + \ldots + p_{11}^2+q^2, \\ p_{i}^2 \in S; \, q^2 \text{ is a prime square.} }} \log p_{1}  \log p_{2}\ldots \log p_{11} \, 2q \, \log q \, \beta(q^2) \label{e3} \\
& \leq 5/4 \sqrt{N} \int_{0}^{1}  \widehat{S}(t)^6 \widehat{S}(-t)^{5}{\psi}(-t)\, dt,   \label{e4} 
\end{align} 

\noindent
in above inequalities: from \eqref{e2} to \eqref{e3} we use the lower bound $ 2/5$ on $ \beta(t^2) $ in the interval $  (N,4N] $ and from \eqref{e3} to \eqref{e4} we use the orthogonality of the functions $ t \mapsto e(nt) $ on $ [0,1]$.

\vspace{2mm}
\noindent
We apply the circle method to estimate the integral on the right hand side of \eqref{32}. To this end,  we set $Q = ( \log N )^{B}  \, A^{48}$, $M = \frac{N}{ (\log N)^{2B}}$, where $B$ is a large absolute constant and, for any integers $a$ and $q$ satisfying  

\vspace{-2mm}
\begin{equation}
\label{maj}
0 \leq a \leq  q \leq Q \;\;\text{and}\;\; (a,q) =1 , 
\end{equation}

\noindent
we call the interval  $[\frac{a}{q} - \frac{1}{M}, \frac{a}{q}+\frac{1}{M})$ the major arc ${\mathfrak M}(\frac{a}{q})$. It is easily checked that distinct major arcs are in fact disjoint when $M > 2Q^2$, which holds when $N$ is sufficiently large depending only on $A$.
We denote by  ${\mathfrak M}$ the union of the family of 
major arcs  ${\mathfrak M}(\frac{a}{q})$. Each interval in the complement of  ${\mathfrak M}$ in $[0, 1)$ is called a minor arc. We denote the union of the minor arcs by ${\mathfrak m}$.

\vspace{2mm}
\noindent
We have 

\vspace{-2mm}
\begin{equation}
\label{32ff}
\int_{0}^{1}  \widehat{S}(t)^6 \widehat{S}(-t)^{5}\psi(-t)\, dt \; = \; 
\int_{-\frac{1}{M}}^{1-\frac{1}{M}}  \widehat{S}(t)^6 \widehat{S}(-t)^{5}\psi(-t)\, dt \;, 
\end{equation} 

\noindent
by the periodicity of the integrand. From the definitions given above it is easily seen  that the interval $[-\frac{1}{M}, 1-\frac{1}{M})$ is the union of ${\mathfrak m}$ and ${\mathfrak M}\setminus [1-\frac{1}{M}, 1+\frac{1}{M})$. Since distinct major arcs are disjoint, it then follows that the right hand side of \eqref{32ff} is the same as  

\vspace{-2mm}
\begin{equation}
\label{32f}
\sum_{1 \leq q \leq Q} \sum_{\substack{0 \leq a < q,\\(a,q) =1. }} \int_{{\mathfrak M}(\frac{a}{q})}  \widehat{S}(t)^6 \widehat{S}(-t)^{5}\psi(-t)\, dt   +
\int_{{\mathfrak m}}  \widehat{S}(t)^6 \widehat{S}(-t)^{5}\psi(-t)\, dt.
\end{equation} 

\noindent
We shall presently estimate each of the two terms in \eqref{32f}. We begin by observing that 

\begin{equation}
\label{esn2}
\int_{0}^{1} |\widehat{S}(t)|^{11} \; dt \ll \, |S|^{9} (\log N )^9 A^3 \; .
\end{equation}

\noindent
In effect, the integral in \eqref{esn2} does not exceed $|S| \, \log N \, E_{5}(S)$. Thus \eqref{esn2} follows from
$|S| \geq N^{\frac{1}{2}}/A  \,(\log N)$ and 

\vspace{-2mm}
\noindent
\begin{align}
E_{5} (S) &= \sum_{\substack{ p_{1}^2 + p_{2}^2 + \ldots + p_{5}^2 = p_{7}^2 + p_{8}^2 + \ldots + p_{10}^2, \\ p_{i}^2 \in S.}} \log p_{1}  \log p_{2}\ldots \log p_{10}  \\
& \leq (\log N)^{10}  \, \sum_{1 \leq n \leq 20N } R_{5}^{2}(n) \ll N^{\frac{3}{2}} (\log N)^{5} |S|^5, \label{esn3}
\end{align}

\noindent
where $R_{5}(n)$ denotes the number of representations of an integer $n$ as a sum of five elements of $S$. 
To verify \eqref{esn3} we note that $ R_{5}(n) = 0$ when $n > 20 N$ and $ R_{5}(n) \leq r_{5}(n)$, the number of representations of $n$ as a sum of five squares of prime numbers, and we have that $r_{5}(n) \ll n^{\frac{3}{2}} / (\log n)^5$  \cite[Theorem 11 ]{hua}, by an application of the circle method. As a consequence of \eqref{esn2} we have 

\vspace{-3mm}
\begin{equation}
\label{dism}
\sum_{1 \leq q \leq Q} \sum_{\substack{0  \leq a < q,\\(a,q) =1. }} \int_{{\mathfrak M}(\frac{a}{q})} |\widehat{S}(t)|^{11} dt \; \leq \; \int_{-\frac{1}{M}}^{1-\frac{1}{M}}  |\widehat{S}(t)|^{11} dt  \ll
|S|^{9}  \, (\log N)^9 \, A^3.
\end{equation}

\subsection{The Minor Arc Contribution}\label{mincon}
Here we bound the second term in \eqref{32f}. Let us first verify that for all $t \in {\mathfrak m}$ we have  

\vspace{-2mm}
\begin{equation}
\label{min0}
|\psi(t)| \ll \frac{N}{A^6} \; ,
\end{equation}

\noindent
when $N$ is large enough, depending only on $A$. Indeed, for any real $t$ Dirichlet's approximation theorem gives a rational number $\frac{a}{q}$ satisfying $|t-\frac{a}{q}| \leq \frac{1}{qM}$ together with $1 \leq q \leq M$ and $(a,q)=1$.  When $t$ is in ${\mathfrak m}$ we see that $\frac{a}{q}$ is in $[0,1]$ since ${\mathfrak m} \subseteq [\frac{1}{M}, 1-\frac{1}{M})$. Consequently, we also have $0 \leq a \leq q$. Since, however, $t$ is not in ${\mathfrak M}$, we must then
have $Q < q$ on account (\ref{maj}). We then conclude using $q^2 \leq qM$ that for each  $t$ in ${\mathfrak m}$ there are integers
$a$ and $q \neq 0$ with $(a,q) =1$ satisfying

\vspace{-2mm}
\begin{equation}
\label{dir}
|t- \frac{a}{q}|\leq \frac{1}{q^2} \;\;\text{ and }\;\; Q < q \leq M.
\end{equation}

\noindent
To get a bound on $ \psi(t)$ when $ t \in {\mathfrak{m}}$, we appeal to the following lemma.

\begin{lem}
\label{exposum}
Let $ \alpha $ be a real number such that 
\begin{equation*}
\alpha = \frac{a}{b} + \lambda, \, \,  (a,q)=1, \, \,  |\lambda| \leq \frac{1}{q^2}, \, \, Q < q \leq M.
\end{equation*}

\noindent
and let $ T(u) = \sum_{0 \leq n \leq u}  \mathbbm{1}_{P}(n)  \, \log{n} \, e(n^2 \alpha)$. Then we have

\begin{equation}
\label{tumax}
\max_{0 \leq u \leq \sqrt{5N}} \, | T(u)| \ll \, \frac{\sqrt{N}}{A^6} \, .
\end{equation}
\end{lem}

\vspace{2mm}
\noindent 
{\sc Proof.---} On the assumption on $ \alpha $, we have 

\begin{equation}
\label{tubound}
\sum_{x < n \leq 2x}  \mathbbm{1}_{P}(n)  \, \log{n} \, e(n^2 \alpha) \ll \, x  \, (\log{x})^{c} \left( \frac{1}{Q} + \frac{1}{x^{1/2}} + \frac{M}{x^2} \right)^{\frac{1}{8}},
\end{equation}

\noindent
for some absolute constant $ c > 0 $; see \cite[Lemma 2.1]{chen}, for example. From this it follows that 

\begin{equation}
\label{tub}
 T(x) = \sum_{0 \leq n \leq x}  \mathbbm{1}_{P}(n)  \, \log{n} \, e(n^2 \alpha) \ll \, x  \, (\log{x})^{c+1} \left( \frac{1}{Q} + \frac{1}{x^{1/2}} + \frac{M}{x^2} \right)^{\frac{1}{8}},
\end{equation}

\noindent
by dividing the interval $ [0,x] $ into dyadic intervals $ (\frac{x}{2^{j+1}}, \frac{x}{2^{j}}]$;  $ j =  0,1,\ldots,\log{x}$  and using the fact that the right of \eqref{tubound} is increasing function of $x$.

\vspace{2mm}
\noindent
Again using the fact that the right of \eqref{tub} is increasing function and recalling the values of $ Q, M $; we get that

\begin{equation}
\max_{0 \leq u \leq \sqrt{5N}}|T(u)|  \ll \, \sqrt{N}  \, (\log{N})^{c+1} \left( \frac{1}{(\log{N})^{B} A^{48}} + \frac{1}{N^{1/4}} + \frac{1}{(\log{N})^{2B}} \right)^{\frac{1}{8}}.
\end{equation}
 
\vspace{2mm} 
\noindent
For large values of $N$ depends on $A$, large absolute value of 
$B$ depends on $c$, we then get

\begin{equation}
\max_{0 \leq u \leq \sqrt{5N}}|T(u)|  \ll \frac{\sqrt{N}}{A^6},
\end{equation}

\noindent
this proves the lemma.

\qed

\vspace{2mm}
\noindent
Now we return to bound $\psi(t)$ on minor the arcs. By the Properties of Riemann-Stieltjes integral we have

\begin{equation}
\psi(t) = \int_{0}^{\sqrt{5N}} \, \, 2 u \, \beta(u^2) \, d T(u),
\end{equation}

\noindent
where $T(u)$ defined as in Lemma \ref{exposum}. Thus, on integrating by parts and using the inequality \eqref{tumax}, we have

\begin{equation}
| \psi(t) \ll \sqrt{N} \max_{0 \leq u \leq \sqrt{5N}}|T(u)| \ll \frac{N}{A^6},
\end{equation}

\noindent
on remarking that $ 2 \, u \, \beta(u^2) $ is piecewise monotonic in the interval $ [0, \sqrt{5N}]$.

\vspace{2mm}
\noindent
From \eqref{min0} and \eqref{esn2} it now follows that for all $N$ large enough, depending only on $A$, we have 

\vspace{-2mm}
\begin{equation}
\label{min1}
\int_{{\mathfrak m}} |\widehat{S}(t)|^{11}\,|\psi(t)| \,  \; dt \ll \frac{{N}}{A^{6}} \int_{0}^{1} |\widehat{S}(t)|^{11} \; dt \;\ll \frac{N|S|^9 (\log {N})^9}{A^3} \ll \frac{|S|^{11} (\log{N})^{11}}{A} ,
\end{equation}

\noindent
since $|S| \geq N^{\frac{1}{2}}/A \, \log{N}$. An application of the triangle inequality now allows us to conclude that 

\vspace{-2mm}
\begin{equation}
\label{min}
\int_{{\mathfrak m}}  \widehat{S}(t)^6 \widehat{S}(-t)^{5}\psi(-t)\, dt \;\ll\; \frac{|S|^{11} (\log{N})^{11}}{A}.
\end{equation}

\subsection{The Function $\psi$ on a Major Arc}

\noindent
Let us set $W = 2U$, where $U$ is as defined staring of Section \ref{seccir}. For any integers $a$, $q$ and $r$, with $ q >0$, we set $V_{q}(a,r) = \sum_{ \substack{0 \leq m < q, \\(r+mW,qW)=1} } e\left(\frac{a(r+mW)^2}{q}\right)$.

\begin{prop} Let $a$ and $q$ be any integers satisfying \eqref{maj}. Then for all $t$ in the major arc  ${\mathfrak M}(\frac{a}{q})$ we have  

\vspace{-2mm}
\begin{equation}
\label{maj1}
\psi(t) = \frac{1}{\phi(qW)} \sum_{\substack{0 \leq r < W, \\(r, W) = 1.}} V_{q}(a,r) \; \overline{\widehat{\beta}}\left(t-\frac{a}{q}\right) + O\left(\phi(W)N \exp{(-c \sqrt{\log N})}\right).
\end{equation}

\end{prop}

\vspace{2mm}
\noindent
{\sc Proof.---} Let $\theta = t-\frac{a}{q}$ and  $f(u) = 2u \log{u}\beta(u^2)e(u^2\theta)$ for any real $u$. we have

\vspace{-2mm}
\begin{equation}
\label{maj3}
\psi(t) = \sum_{\substack{0 \leq r < W, \\(r, W) = 1.}}\, 
 \,\sum_{n \, \equiv \, r\, {\rm mod}\, W} \mathbbm{1}_{P}(n) \, n \, 2n \, \log{n}  \beta(n^2) e(n^2t)  + O \left( 9^{A^{l}}\right), 
\end{equation}

\noindent
on recalling expression of $\psi(t)$, noticing the fact that all primes more than $U$ are co-prime to $W$, and trivially estimating contribution to the sum over the interval $[0,U]$ using an upper bound $3^{A^{l}}$ on $U$.

\vspace{2mm}
\noindent
We now find an asymptotic formula for the inner sum on right of  \eqref{maj3}. To this end, we let $X=[\sqrt{5N}]$, and for  $ n \in [1,X] $ let 

\vspace{-2mm}
\begin{equation}
\label{sn}
S_{n} = \sum_{\substack{m \leq n, \\ m \, \equiv \, r {\rm mod} \, W}} \mathbbm{1}_{P}(n) \, e(\frac{a m^2}{q}) \, .
\end{equation} 

\vspace{2mm}
\noindent
Using the fact that $ q \leq Q$, we get 

\begin{equation}
S_{n} = \sum_{\substack{0 \leq m < q, \\ (r+mW,qW)=1}} \, e \left( \frac{a (r+mW)^{2}}{q} \right)  
\sum_{\substack{d, \\ (r+mW) + d qW \leq 1, \\ (r+mW) + d qW \, \, \text{is a prime}}}  1  + \, O (Q) \;.
\end{equation}

\vspace{2mm}
\noindent
As $ n \leq \sqrt{5N} $ and $ qW \leq (\log N)^{B+1}$ for large values of $N$, by appealing to the Siegel–Walfisz theorem, we get the asymptotic expression

\begin{equation}
\label{asymsn}
S_{n} = \frac{Li(n)}{\phi(qW)} \sum_{\substack{0 \leq m < q, \\ (r+mW,qW)=1}} \, \, e \left( \frac{a (r+mW)^{2}}{q} \right)  + \,\, O  
\left( \sqrt{N}  \, \exp{(-c \sqrt{\log N})} \right) \;,
\end{equation}

\noindent
where $Li(n) = \int_{2} ^{n} \frac{dt}{\log{t}} $, and $c$ is a positive absolute constant.

\vspace{2mm}
\noindent
Using the functions $f(u)$ and $S_{n}$, we have
\begin{align}
\sum_{n \, \equiv \, r\, {\rm mod}\, W} \mathbbm{1}_{P}(n) \, n \, 2n \, \log{n}  \beta(n^2) e(n^2t)  &= \sum_{n=2}^{X} \, (S_{n}-S_{n-1}) \, f(n) \label{mean1} \\
                            &=  S_{X} \, f(X+1) + \sum_{n=2}^{X} S_{n} \, (f(n)-f(n+1)) \label{mean2}\; .
\end{align}

\vspace{2mm}
\noindent
As $ |\theta| \leq \frac{(\log{N})^{2B}} {N} $, the Mean-Value theorem implies that 

\begin{equation}
f(n) - f(n+1) \ll (\log{N})^{2B + 1} \;.
\end{equation}

\vspace{2mm}
\noindent
Hence the sum on left of \eqref{mean1} becomes

\begin{equation}
\label{mean3}
\frac{V_{q}(a,r)}{\phi(qW)} \,\, \left[ Li(X) f(X+1) + \sum_{n=2}^{X} Li(n) \, (f(n)-f(n+1)) \right] + O \left( N \exp{(-c \sqrt{\log{N}})} \right) \;.
\end{equation}

\vspace{2mm}
\noindent
As $Li(2)=0$, we then rewrite \eqref{mean3} as

\begin{equation}
\frac{V_{q}(a,r)}{\phi(qW)} \,\, \sum_{n=3}^{X} \int_{n-1}^{n} \frac{f(n)}{\log{x}} \, dx + O \left( N \exp{(-c \sqrt{\log{N}})} \right) \;.
\end{equation}

\vspace{2mm}
\noindent
When $ n-1 < x < n $, the Mean-value theorem reveals that 
\begin{equation}
f(n) = f(x) + O \left( (\log{N})^{2B+1} \right)\;,
\end{equation}

\noindent
and so 
\begin{equation}
\label{mean4}
\sum_{n=3}^{X} \int_{n-1}^{n} \frac{f(n)}{\log{x}} \, dx = \, \int_{2}^{X} 2x \, \beta(x^2) \, e(\theta x^2) dx  + O \left( \sqrt{N} \, (\log{N})^{2B+1}\right) \;.
\end{equation}

\vspace{2mm}
\noindent
Note that the integral on right of above is nothing but  ${\overline{\widehat{\beta}}}(\theta)$. Thus, we have an asymptotic formula for the inner sum on right of \eqref{maj3} as follows

\begin{equation}
\sum_{n \, \equiv \, r\, {\rm mod}\, W} \mathbbm{1}_{P}(n) \, n \, 2n \, \log{n}  \beta(n^2) e(n^2t) = \frac{V_{q}(a,r)}{\phi(qW)} \, {\overline{\widehat{\beta}}}(\theta) + O \left( N \exp{(-c \sqrt{\log{N}})} \right) \;,
\end{equation}

\noindent
for large enough $N$ depends on only on $A$. Substituting this into \eqref{maj3} we get an asymptotic formula for $ \psi(t)$ as in \eqref{maj1}.

\vspace{2mm}
\noindent
We need the following proposition, which provides information about $V_{q}(a,r)$.

\begin{prop} \label{lemg} Let $a$ and $q$ be integers satisfying \eqref{maj} and $r$ any integer co-prime to $W$. Then we have 

\vspace{1mm}
\noindent
$(i)$ $V_{q}(a,r) =0$ unless $q|2W$ or there is a prime $p >w $ such that $p|q$.

\noindent
$(ii)$ $\frac{1}{\phi(qW)}|V_{q}(a,r)| \ll \frac{1}{ \phi(W) w^{97/200}}$ when $q$ does not divide $2W$.

\end{prop}

\vspace{2mm}
\noindent
We prove this Proposition with the help of following lemma.

\begin{lem}
\label{factlem}
Let $ P(z) = c_{0} z^2 + c_{1} z + c_{2}$ be a polynomial with integer coefficients and let $d$ be a positive integer with $d=d_{1} d_{2}$ and $d_{2}$ divides $c_{0}$. Then 

\begin{equation}
\label{fact}
\sum_{0 \leq m < d} e\left(\frac{P(m)}{d}\right) =
 \sum_{0 \leq m_1 < d_1} e\left(\frac{P(m_1)}{d}\right)
\sum_{0\leq m _2 < d_2} e\left(\frac{c_1m_2}{d_2}\right).
\end{equation}

\end{lem}

\vspace{2mm}
\noindent
{\sc Proof.---} See \cite[page 26]{gsr}, for example. 

\vspace{2mm}
\noindent
We now give a proof of the above Proposition with the aid of this lemma.
Since $(r,W)=1$, the condition $(r+mW,qW)=1$ in the definition of $V_{q}(a,r)$ can be replaced by the condition $(r+mW,q)=1$, thus we have

\begin{equation}
\label{vqr1}
V_{q}(a,r) = \sum_{ \substack{0 \leq m < q, \\(r+mW,q)=1} } e\left(\frac{a(r+mW)^2}{q}\right) \;.
\end{equation}

\vspace{2mm}
\noindent
Let $q=UV$, where $U$ is $w$-smooth and $(V,W)=1$ and let $a=a_{2}U+a_{1}V, (a_{1},U)=1$ and $ (a_{2},V)=1$. Then from \eqref{vqr1} follows that

\begin{equation}
\label{vqr2}
e(\frac{-a r^2}{q}) \,  V_{q}(a,r) = \left( \sum_{ \substack{0 \leq m_{2} < V, \\(r+m_{2} W,V)=1} } e\left(\frac{a_{2} W^2 m_{2}^2 + 2 a_{2} W r m_{2}}{V}\right) \right) \, \left( \sum_{ \substack{0 \leq m_{1} < U, \\(r+m_{1} W,U)=1} } e\left(\frac{a_{1} W^2 m_{1}^2 + 2 a_{1} W r m_{1}}{U}\right) \right) \;.
\end{equation}

\vspace{2mm}
\noindent
Now we analyze the second term in the product of right of the above equation. Since $U$ is $w$-smooth and $(r,W)=1$, the condition $(r+m_{1} W, U)=1$ is always holds, thus we get

\begin{equation}
\label{vqr3}
\sum_{ \substack{0 \leq m < U, \\(r+mW,U)=1} } e\left(\frac{a_{1} W^2 m^2 + 2 a_{1} W r m}{U} \right) = \sum_{0 \leq m < U} e\left(\frac{a_{1} W^2 m^2 + 2 a_{1} W r m}{U} \right) \,.
\end{equation}

\vspace{2mm}
\noindent
We write $U=U_{1} U_{2}$, where $U_{1} = \frac{U}{(U,W^2)},\, U_{2}=(U,W^2)$. Note that $U_{2}| \,a_{1} W^2$, thus applying the Lemma \ref{factlem}, we get that

\begin{equation}
\label{vqr4}
\sum_{0 \leq m < U} e\left(\frac{a_{1} W^2 m^2 + 2 a_{1} W r m}{U} \right) = \left( \sum_{0 \leq m_{1} < U_{1}} e\left(\frac{a_{1} W^2 m_{1}^2 + 2 a_{1} W r m_{1}}{U} \right) \right) \,  \left( \sum_{0 \leq m_{2} < U_{2}} e\left(\frac{2 a_{1} W r m_{2}}{U_{2}} \right) \right) \,.
\end{equation}

\vspace{2mm}
\noindent
We can conclude from \eqref{vqr4},\eqref{vqr3} and \eqref{vqr2} that $V_{q}(a,r)=0$ unless $U_{2}| 2a_{1}Wr$. That is, unless $(U,W^2)|2a_{1}Wr$ we have$V_{q}(a,r)=0$.

\vspace{2mm}
\noindent
Since $(a_{1}, U)=1$ and $(r,W)=1$, it follows that $V_{q}(a,r)=0$ unless $(U,W^2)|2W$. We note that $(U,W^2) = (q,W^2)$ and that $(q,W^2)|2W$ is equivalent to $\inf (v_{p}(q), 2v_{p}(W)) \leq v_{p}(2W)$ for all primes $p|2W$. From the definition of $W$ we have $2v_{p}(W) > v_{p}(2W)$ for all primes $p| 2W$. Consequently, $V_{q}(a,r) = 0$ unless $v_{p}(q) \leq v_{p}(2W)$ for all primes $p|2W$, which is the same as $(i)$.

\vspace{2mm}
\noindent
To prove $(ii)$, we may assume that $(q,W^2) | 2W$ and $(q,W^2)|2W$ and $V>1$. We can conclude from $(q,W^2)=(U,W^2)$ and $(q,W^2)|2W$ that $U|2W$, inparticular $U|W^2$. Thus, we have $(U,W^2)=U$. It follows that

\begin{equation}
\label{vqr5}
\sum_{ \substack{0 \leq m_{1} < U, \\(r+m_{1} W,U)=1} } e\left(\frac{a_{1} W^2 m_{1}^2 + 2 a_{1} W r m_{1}}{U}\right) = U \;,
\end{equation}

\noindent
again using the same fact that $(r+m_{1}W,U)=1$ is always holds, as $U$ is a $w$-smooth number. On combining \eqref{vqr2} and \eqref{vqr5} we get that

\begin{equation}
\label{vqr6}
e(\frac{-a r^2}{q}) \,  V_{q}(a,r) =  U  \sum_{ \substack{0 \leq m_{2} < V, \\(r+m_{2} W,V)=1} } e\left(\frac{a_{2} W^2 m_{2}^2 + 2 a_{2} W r m_{2}}{V}\right) \;.
\end{equation}

\vspace{2mm}
\noindent
Multiplying by a function $ e(\frac{-a_{2} r^2}{V}) $ both side of above equation and change the variable $ r+m_{2}W \mapsto x $ in the summation on right of \eqref{vqr6} and using the fact that $(V,W)=1$, we get 

\begin{equation}
\label{vqr7}
e\left((-r^2 (a_{2}/V + a/q)\right) \,  V_{q}(a,r) =  U  \sum_{ \substack{0 \leq x < V, \\(x,V)=1} } e\left(\frac{a_{2} x^2}{V}\right) \,. 
\end{equation} 

\vspace{2mm}
\noindent
We have a following bound on the summation on right of \eqref{vqr7}

\begin{equation}
\label{vqr8}
\sum_{ \substack{0 \leq x < V, \\(x,V)=1} } e\left(\frac{a_{2} x^2}{V}\right) \ll_{\epsilon} V^{\frac{1}{2} + \epsilon},
\end{equation}

\noindent
see \cite[Lemma 8.5]{hua}, for example.

\vspace{2mm}
\noindent
From \eqref{vqr8} and \eqref{vqr7}, it follows that 

\begin{equation}
\label{vqr9}
\frac{1}{\phi(qW)}|V_{q}(a,r)| \ll_{\epsilon} \frac{U \,V^{\frac{1}{2}+\epsilon}}{ \phi(qW)} =  \frac{U \,V^{\frac{1}{2}+\epsilon}}{ \phi(V)\phi(UW)} \,,
\end{equation} 

\noindent
here we use the fact that $(V,W)=1$ in the equality on the right of the above equation. Since $U|2W$, we have $\phi(UW) = U \phi(W)$ and we have the lower bound on $\phi(V)$, namely $\phi(V) \gg V / \log log V$. From this and \eqref{vqr9} we get 

\begin{equation}
\label{vqr10}
\frac{1}{\phi(qW)}|V_{q}(a,r)| \ll_{\epsilon} \frac{ \log \log V}{ \phi(W) V^{\frac{1}{2} - \epsilon}} \,.
\end{equation} 

\vspace{2mm}
\noindent
Taking $\epsilon = 1/200$ and using the bound  $\log \log V \leq V^{1/100}$ in \eqref{vqr10} we conclude $(ii)$.

\subsection{The Major Arc Contribution}
\label{majcon}

\noindent
In this subsection we reduce the problem of bounding $E_{6}(S)$ to a finite problem. Let us first dispose of the first term in \eqref{32f}, which we denote here by $T$. Then on writing $T_1$ for
  
\vspace{-2mm}
\begin{equation}
\label{maj11}
\sum_{\substack{0 \leq r < W, \\(r, W) = 1.}}  \sum_{1 \leq q \leq Q} \frac{1}{\phi(qW)} \sum_{\substack{0 \leq a < q,\\(a,q) =1. }}
{V_{q}}(-a,r) \;\int_{{\mathfrak M}(\frac{a}{q})} {\widehat{\beta}}\left(t-\frac{a}{q}\right) \widehat{S}(t)^6 \widehat{S}(-t)^{5}\, dt 
\end{equation} 

\noindent
we deduce by substituting the complex conjugate of right hand side of \eqref{maj1} for $\psi(-t) = \overline{\psi}(t)$ in $T$ and using the triangle inequality together with \eqref{dism} that  

\vspace{-2mm}
\begin{equation}
\label{maj111}
T-T_{1} \ll \phi(W) \,  N \exp{(-c \, \sqrt{\log N})} \int_{0}^{1} |\widehat{S}(t)|^{11} dt \ll  \phi(W) \,  N \exp{(-c \, \sqrt{\log N})} 
\, |S|^9  \, (\log N)^9  \, A^3 \,.
\end{equation} 

\noindent
If we now set 

\vspace{-2mm}
\begin{equation}
\label{maj12}
T(W) = \sum_{\substack{0 \leq r < W, \\(r, W) = 1.}}  \sum_{ q|2W} \frac{1}{\phi(qW)}  \sum_{\substack{0 \leq a < q,\\(a,q) =1. }}
V_{q}(-a,r) \;\int_{{\mathfrak M}(\frac{a}{q})} {\widehat{\beta}}\left(t-\frac{a}{q}\right) \widehat{S}(t)^6 \widehat{S}(-t)^{5}\, dt .
\end{equation}

\noindent
then by $(ii)$ of Lemma \ref{lemg} combined with the triangle inequality and \eqref{dism} we get

\vspace{-2mm}
\begin{equation}
\label{maj13}
T_{1} - T(W) \; \ll \;  \frac{\|\widehat{\beta}\|_{\infty} |S|^{9} 
(\log N)^9 A^3}{w^{97/200}} \ll \frac{A^3 |S|^9 (\log N)^9 N}{w^{97/200}},
\end{equation}

\noindent
 since $\|\widehat{\beta}\|_{\infty} = \sup_{t \in {\bf R}}|\widehat{\beta}(t)|\leq \frac{5N}{2}$. From \eqref{maj13}, \eqref{maj111} and on recalling that $|S| \geq \frac{\sqrt{N}}{A \log N}$ and $w = A^{25} $ we conclude that 

\vspace{-2mm}
\begin{equation}
\label{maj14}
T =  T(W) + O\left( \frac{|S|^{11}}{A}\right),
\end{equation} 

\noindent
when $N$ is sufficiently large, depending only on $A$. Let us now estimate $T(W)$. When $q|2W$ we have $ \phi(qW) = q \phi(W)$ and $(r+mW)^2 \equiv r^2$ modulo $q$  for all integers $m$ and the condition 
$(r+mW)=1 $ holds always. Therefore we have $V_{q}(a,r) = qe\left(-\frac{ar^2}{q}\right)$ when $q|2W$, for all $0 \leq a <q$. Furthermore, since $r \mapsto r+W$ is a bijection from the integers co-prime to $2W$  in $[0,W)$ to those in $(W, 2W]$ co-prime to $2W$,  we obtain  

\vspace{-2mm}
\begin{equation}
\label{co2w}
\frac{1}{\phi(qW)} \sum_{\substack{0 \leq r < W, \\(r, W) = 1.}} V_{q}(-a,r)  = \frac{1}{2 \phi(W)} \sum_{\substack{0 \leq r < 2W, \\(r, 2W) = 1.}} e\left(-\frac{ar^2}{q}\right)
\end{equation}

\noindent
for any $q|2W$ and all $0 \leq a <q$. Also, we have $\widehat{S}(t)^6 \widehat{S}(-t)^{5} = \sum_{x \in S^{11}}  \log{x_{1}} \log{x_{2}} \ldots \log{x_{11}} \, e\left(f(x)t\right)$, where $f(x)$ denotes $x_{1}+x_{2}+ \ldots + x_{6}-x_{7}-\ldots - x_{11}$ for any $x=(x_{1},\ldots,x_{11}) \in S^{11}$. By means of the change of variable $t -\frac{a}{q} \mapsto t$ in the integrals in \eqref{maj12} we then see that 

\vspace{-2mm}
\begin{equation}
\label{maj15}
T(W) = \frac{1}{2 \phi(W)} \sum_{\substack{0 \leq r < 2W, \\(r, 2W) = 1.}} \sum_{ q|2W} \sum_{\substack{0 \leq a < q,\\(a,q) =1. }} \int_{-\frac{1}{M}}^{\frac{1}{M}}  {\widehat{\beta}}\left(t\right) \sum_{x \in S^{11}} \prod_{i=1}^{11} \log{x_{i}}  \, e\left(tf(x)\right) \, e\left(\frac{a(f(x) -r^2)}{q} \right)\, dt .
\end{equation}

\noindent
Finally, on interchanging summations and remarking that

\vspace{-2mm}
\begin{equation}
\label{maj16}
 \frac{1}{2W} \sum_{ q|2W} \sum_{\substack{0 \leq a < q,\\(a,q) =1. }}
e\left(\frac{a(f(x) -r^2)}{q} \right) =  \frac{1}{2W} \sum_{{0 \leq a<2W}}e\left(\frac{a(f(x) -r^2)}{2W} \right)
\end{equation}

\noindent
we conclude that the right hand side of \eqref{maj15} is the same as the left hand side of    

\vspace{-2mm}

%example of equation spliting
%\begin{dmath}
%\label{maj17}
%\frac{W}{\phi(W)}\sum_{\substack{0 \leq r < 2W, \\(r, 2W) = 1.}} \sum_{\substack{x \in S^{11}, \\ f(x)\equiv r^2 \rm{mod}\, 2W}} \log{x_{1}} \log{x_{2}} \ldots \log{x_{11}}  \int_{-\frac{1}{M}}^{\frac{1}{M}}  {\widehat{\beta}}\left(t\right) e\left(tf(x)\right) dt\; \leq  \; \frac{W (\log N)^{11}}{\phi(W)}
%\sum_{\substack{0 \leq r < 2W, \\(r, 2W) = 1.}} \sum_{\substack{x \in S^{11}, \\ f(x)\equiv r^2 \rm{mod}\, 2W}} 1,
%\end{dmath}
\begin{equation}
\begin{split}
\label{maj17}
\frac{W}{\phi(W)}\sum_{\substack{0 \leq r < 2W, \\(r, 2W) = 1.}} \sum_{\substack{x \in S^{11}, \\ f(x)\equiv r^2 \rm{mod}\, 2W}} \log{x_{1}} \log{x_{2}} \ldots \log{x_{11}}  \int_{-\frac{1}{M}}^{\frac{1}{M}}  {\widehat{\beta}}\left(t\right) e\left(tf(x)\right) dt\;\\ \leq  \; \frac{W (\log N)^{11}}{\phi(W)}
\sum_{\substack{0 \leq r < 2W, \\(r, 2W) = 1.}} \sum_{\substack{x \in S^{11}, \\ f(x)\equiv r^2 \rm{mod}\, 2W}} 1,
\end{split}
\end{equation}

\noindent
where we have used $|\int_{-\frac{1}{M}}^{\frac{1}{M}}  {\widehat{\beta}}\left(t\right) e\left(tf(x)\right) dt| \leq \int_{\bf R} \widehat{\alpha}(t) dt = 1$, since
$|{\widehat{\beta}}\left(t\right)| = \widehat{\alpha}(t)$ for all $t \in {\bf R}$. For each invertible square $b$ in ${\bf Z}/{2W}{\bf Z}$, the number of $r$ in $[0, 2W)$ co-prime  to $2W$ and such that $r^2 \equiv b$ modulo $2W$ is $2\tau(U)$. Then it follows from \eqref{maj17} and \eqref{maj15} that 

\vspace{-2mm}
\begin{equation}
\label{maj18}
T(W) \leq  \frac{2 W \tau(U) (\log N)^{11}}{\phi(W)} \left|\{ x \in S^{11}\,|\, f(x) \;\text{an invertible square mod $2W$} \} \right|. 
\end{equation}   

\noindent
On combining \eqref{maj18} with \eqref{maj14}, \eqref{min} and recalling that \eqref{32f} is the same as the integral in \eqref{32}, we get

\begin{equation}
\label{finitepr}
E_{6}(S) \ll \frac{2 W \tau(U) (\log N)^{11}}{\phi(W)} \left|\{ x \in S^{11}\,|\, f(x) \;\text{an invertible square mod $2W$} \} \right|  \,
+ O \left( \frac{|S|^{11} (\log N)^{11}}{A} \right)
\end{equation}

\subsection{Proof of Theorem \ref{main2} Completed}\label{maincom}

\noindent
It remains only to bound the cardinality of the set $\{ x \in S^{11}\,|\, f(x) \;\text{an invertible square mod $2W$} \}$. We find an upper bound the cardinality of this set using Theorem \ref{th4}. Let ${\mathcal Z}$ be the set of integers $n >0$ such that $n^2 \in S$. The set ${\mathcal Z}$ is contained in $[\sqrt{N}, 2 \sqrt{N})$  and satisfies $|{\mathcal Z}|  \geq \frac{\sqrt{N}}{A \log N}$ and $|\{ z \in {\mathcal Z}| z \equiv a \,{\rm mod}\, U\}| \leq \frac{3 \sqrt{N}}{\phi(U) \log N}$, when $N$ is sufficiently large depending on $A$. Finally, let $I = S^{9}$ and for any $x = (x_1, x_2, \ldots, x_9) \in S^{9}$ we set $c(x) = x_1 +  \ldots +x_4 -x_5- \ldots -x_9$. Then with ${R}_{U}({\mathcal Z}, {\bf c})$ as in Theorem \ref{th4} we have that 

\vspace{-2mm}
\begin{equation}
\label{maj19}
 \left|\{ x \in S^{11}\,|\, f(x) \;\text{an invertible square modulo $2W$} \} \right| \; \leq \; |{R}_{U}({\mathcal Z}, {\bf c})|,
\end{equation}

\noindent
since $U |2W$.
On combining the bound for $|{R}_{U}({\mathcal Z}, {\bf c})|$ given by Theorem \ref{th4} with \eqref{maj19} and \eqref{finitepr} and after noticing that $\frac{W}{\phi(W)} (\frac{U}{\phi(U)})^2 \ll (\log A)^3$ we finally obtain \eqref{three11}, as required.

\section{Monochromatic Representation}
\label{fin}

\noindent
Here we deduce Theorem \ref{th1} from Theorem \ref{main2}. Before this deduction we give a lemma with the following notation. For any subset $S$ of the integers, we write $e_{6}(S)$ for the number of tuples $(x_{1},x_{2},\ldots,x_{12})$ in $S^{12}$ satisfying $x_{1}+ \ldots + x_{6} = x_{7}+ \ldots + x_{12}$. We observe that if $S \subset [N,4N]$ satisfying the hypothesis of Theorem \ref{main2}, then we can conclude from \eqref{three11} that 

\begin{equation}
\label{adenerbdd}
e_{6}(S) \ll \frac{|S|^{11}}{N^{\frac{1}{2}} \log N} \exp\left(\frac{\left(3\log 2 + o(1)\right)\log A}{\log \log A}\right) \,.
\end{equation}

\noindent
Now we state the lemma as follows.

\begin{lem}
\label{sarkfinite}
Let $N$ be positive integer and let $D \geq 1$ be a real number satisfying the condition $N \geq 72D+12$. If $S$ be a subset of the interval $(N,4N]$ such that
\begin{equation}
\label{e6sbdd}
e_{6}(S) \leq \frac{|S|^{12} D}{3N}
\end{equation}

\noindent
and if $S$ contains an integer that is not divisible by any prime $p \leq \lceil 6D \rceil$ then every integer $n \geq 30N (2 \lceil 6D \rceil +1)$ is a sum of no more than $\frac{n}{N}$ elements of $S$.

\end{lem}

\vspace{2mm}
\noindent
{\sc Proof.---} See \cite[Lemma 1.2]{gsr}, for example.

\vspace{2mm}
\noindent
We now give the proof of Theorem \ref{th1}.
Since $s(K)$ is increasing with $K$, it suffices to prove Theorem \ref{th1} for all $K$ sufficiently large. For such a $K$, 
let $\cup_{1 \leq i \leq K} {\mathfrak Q}_i$ be a partition 
of the set of square of primes ${\mathfrak Q}$ into $K$ disjoint subsets. 

\vspace{2mm}
\noindent
We set $N_{0} = 2 K^{50}$. Let $N$ be an integer $\geq N_{0}$. There is an $i, 1 \leq \leq K$, such that $\mathfrak Q_{i} \cap (N,4N]$ contains atleast $\frac{\sqrt{N}}{K \log N}$ elements of $\mathfrak Q$. For such $i$ we set $S=\mathfrak Q_{i} \cap (N,4N]$. Then $S$ is set of square of primes in $(N,4N]$ with $|S| \geq \frac{\sqrt{N}}{K \log N}$ and no integer in $S$ is divisible by a prime $p \leq K^{25}$. It now follows from \eqref{adenerbdd} that \eqref{e6sbdd} holds with $D \ll K \exp\left(\frac{\left(3\log 2 + o(1)\right)\log K}{\log \log K}\right) $. Since elements of $S$ does not divisible by any prime $p \leq \lceil 6D \rceil$ when $N$ is large enough, we may apply Lemma \ref{sarkfinite} to $S$ to deduce that  every integer $n \geq (288D+72)N$ is a sum of no more than $\frac{n}{N}$ elements of $S$. In particular, there is a $C_1 >0$ such that every integer in $I(N) = ((288D+72)N, (288D +73)N]$ is a sum of at most $C_1 D$ squares of primes all belonging to $S$ and therefore to ${\mathfrak Q}_i$. Thus for all large enough $N$, every integer in the interval $I(N)$ can be expressed as a sum of no more than $C_1D$ squares of primes all of the same colour. On remarking that the interval $I(N)$ meets $I(N+1)$ for all large enough $N$, we obtain that $s(K) \leq C_1 D$. This yields the conclusion of Theorem \ref{th1} since $C_1 D \ll K\exp\left(\frac{\left(3\log 2 + o(1)\right)\log K}{\log \log K}\right)$ .

\vspace{3mm}
\noindent
{\bf Acknowledgement :} This work was carried out under the CEFIPRA project  5401-1.

\vspace{5mm}

\end{document}